\documentclass[10pt,reqno]{amsart}

\usepackage{amssymb, amsthm, amsmath, latexsym}
\usepackage{amsfonts}

\newcommand{\bp}{\begin{proof}}
\newcommand{\ep}{\end{proof}}
\newcommand{\be}{\begin{equation}}
\newcommand{\ee}{\end{equation}}

\newcommand{\e}{\nu}
\newcommand{\de}{\varepsilon}
\newcommand{\R}{\mathbb{R}}

\newtheorem{thm}{Theorem}[section]

\numberwithin{equation}{section}
\newtheorem{remark}{Remark}[section]
\renewenvironment{proof}[1][Proof]{\textbf{#1.} }{\ \rule{0.5em}{0.5em}}

\numberwithin{equation}{section}

\begin{document}
\title[Solutions for a Nonlocal Conservation Law with Fading Memory]{Solutions
 for a Nonlocal Conservation Law with Fading Memory}
\author{Gui-Qiang Chen \and Cleopatra Christoforou}
\address{Department of Mathematics\\
Northwestern University\\
2033 Sheridan Road\\
Evanston, Illinois 60208, USA}

\email{gqchen@math.northwestern.edu}
\urladdr{http://www.math.northwestern.edu/$\sim$gqchen}

\email{cleo@math.northwestern.edu}
\urladdr{http://www.math.northwestern.edu/$\sim$cleo}


\date{March 30, 2006}
%
% Abstract
\begin{abstract}
Global entropy solutions in $BV$ for a scalar nonlocal conservation
law with fading memory are constructed as limits of vanishing
viscosity approximate solutions. The uniqueness and stability of
entropy solutions in $BV$ are established, which also yield the
existence of entropy solutions in $L^\infty$ while the initial data
is only in $L^\infty$. Moreover, if the memory kernel depends on a
relaxation parameter $\de>0$ and tends to a delta measure weakly as
measures when $\de\to 0+$, then the global entropy solution sequence
in $BV$ converges to an admissible solution in $BV$ for the
corresponding local conservation law.
\end{abstract}
\maketitle

\section{Introduction and Main Theorems}

We study global entropy solutions to a scalar nonlocal conservation
law with fading memory: \be \label{1.fad}
u_t+f(u)_x+\int_0^tk(t-\tau)f(u(\tau))_x\,d\tau=0, \hspace{1cm}
x\in\mathbb{R}, \ee and initial data \be \label{1.in-fad} u(0,
x)=u_0(x), \ee where $f:\mathbb{R}\to\mathbb{R}$ is a smooth
function and $u_0\in BV(\mathbb{R})$. For simplicity, we
sometimes use the notation $u(t):=u(t,x)$ to emphasize the state at
$t>0$ as in \eqref{1.fad}.

In one-dimensional viscoelasticity, hyperbolic conservation laws
\be\label{1.cl} U_t+F(U)_x=0 \ee correspond to the constitutive
relations of an elastic medium when the value of the flux function
$F$ at $(t,x)$ is solely determined by the value of $U(t,x)$.
However, this model \eqref{1.cl} is inadequate when viscosity and
relaxation phenomena are present. In that case, the flux function
depends also on the past history of the material, i.e. on
$U(\tau,x)$ for $\tau<t$. Under these circumstances, we say that the
material has memory. An important class of media of this type are
materials with fading memory, which correspond to the constitutive
relations with flux functions of the form \be\label{1. fad mem}
F(U(t,x))+\int_0^t k(t-\tau) G(U(\tau,x))\,d\tau, \ee where $F$, $G$
are smooth functions and $k$ is a smooth kernel, integrable over
$\mathbb{R}_+:=[0,\infty)$. When the kernel satisfies appropriate
conditions motivated by physical considerations, the influence of
the memory term in \eqref{1. fad mem} reflects a damping effect.
Consequently, global smooth solutions exist for given small initial
data (cf. Renardy-Hrusa-Nohel \cite{RHN}), in contrast to the
situation with elastic media in which classical solutions in general
break down in finite time even when the initial data is small.
However, when  the initial data is large, the destabilizing action
of nonlinearity of the flux function $f$ prevails over the damping,
and solutions break down in a finite time; see Dafermos \cite{D2}
and Malek-Madani--Nohel \cite{MN}.

In this paper we first construct global entropy solutions in $BV$ to
the nonlocal conservation law \eqref{1.fad} with fading memory via
the vanishing viscosity approximation by adding the artificial
viscosity term as follows:
\be \label{1.fad-vis}
u^{\e}_t+f(u^{\e})_x+\int_0^tk(t-\tau)f(u^{\e}(\tau))_x\,d\tau =\e\,
(u^{\e}+\int_0^t k(t-\tau) u^{\e}(\tau)\,d\tau)_{xx}. \ee For a
scalar local conservation law modeling elastic materials, such a
result was first established in \cite{K,O,Volpert}.

The main motivation for the vanishing viscosity approximation
\eqref{1.fad-vis} is that the conservation law \eqref{1.fad} can be
viewed as a linear Volterra equation, which was first observed by
MacCamy \cite{MC} and later employed in Dafermos \cite{D3} and
Nohel--Rogers--Tzavaras \cite{NRT} (also see \cite{CD}). In this
way, it is easy to extract the damping character of the memory term.
Let $r(t)$ be the resolvent kernel associated with $k(t)$:
\be\label{1.resolvent relation} r+k* r=-k. \ee Then we can write
\eqref{1.fad} as
\[-f(u)_x=u_t+\int_0^tr(t-\tau)\, u'(\tau)\,d\tau.\]
Integrating by parts yields \be\label{1.res-fad}
u_t+f(u)_x+r(0)u=r(t) u_0-\int_0^t r'(t-\tau) u(\tau)\,d\tau,\ee
which is equivalent to \eqref{1.fad}. The vanishing viscosity
approximation \eqref{1.fad-vis} is equivalent to the following
artificial viscosity approximation to \eqref{1.res-fad}:
\be\label{1.res-fad-vis} u_t^{\e}+f(u^{\e})_x+r(0)u^{\e} =r(t)
u_0-\int_0^t r'(t-\tau) u^{\e}(\tau)\,d\tau+\e\, u^{\e}_{xx}.\ee
Note that the above argument applies only if the nonlinearity $f$ in
the instantaneous response is the same as in the memory term chosen
in \eqref{1.fad}; also see \cite{D3,MC,NRT} for the same
restriction.
The artificial viscosity term in \eqref{1.fad-vis} is chosen so that
\eqref{1.res-fad-vis} has the standard  artificial viscosity term to
ensure the $L^\infty$ and $BV$ estimates of $u^\nu$ (also see
\cite{CD,D3,NRT}).

The existence of a unique, regular local solution $u^\e(t,x)$ of
\eqref{1.res-fad-vis} when the initial data $u_0$ is smooth can be
established through the standard Banach Fixed Point Theorem. The
local solution may be extended to a global solution with the help of
the apriori $L^\infty$ estimate established in Section 2.1.

A function $u=u(t,x)$ is called an {\it entropy solution} to the
Cauchy problem \eqref{1.fad}--\eqref{1.in-fad} if it satisfies that,
for any test function $\varphi\in C_0^1(\mathbb{R}_+^2)$ with
$\mathbb{R}_+^2:=\mathbb{R}_+\times \mathbb{R}$, $\varphi\ge 0$,
\begin{align}
&\iint_{\mathbb{R}^2_+}\big(\eta(u)\varphi_t +q(u)\varphi_x
+\eta'(u)(r(t)u_0-r(0)u-\int_0^t
r'(t-\tau)u(\tau)d\tau)\varphi\big)dt dx\nonumber\\
&\qquad+\int_{\mathbb{R}}\eta(u_0(x))\varphi(0,x)dx\ge 0,
\label{entropy-ineq}
\end{align}
for any convex entropy $\eta(u),$ where $q(u)$ is
the corresponding entropy flux satisfying $q'(u)=\eta'(u)f'(u)$.

Before we state the results, we introduce some notations. Let
$\varrho$ be the standard mollifier. We define the mollification of
$u_0$ to be \be\label{1.in-z} u_0^{\nu}:=(u_0 \chi_\nu)*
\varrho_{\nu},\ee where, for each $\nu>0$,
$\varrho_{\nu}(x):=\frac{1}{\nu}\, \varrho(\frac{x}{\nu})$ and
$\chi_\nu(x):=1$ for $|x|\le 1/\nu$ and $0$ otherwise. The main
result is the following.

\begin{thm}\label{thm1}
Consider the Cauchy problem \eqref{1.fad-vis} with Cauchy data: \be
\label{1.in-fad-vis} u^{\e}(0, x)=u_0^\nu(x),\ee where the initial
data $u_0^\nu$ is given by \eqref{1.in-z} and $u_0\in
BV(\mathbb{R})$. Let the resolvent kernel $r$ associated with $k$ as
defined in \eqref{1.resolvent relation} be a nonnegative,
non-increasing function in $L^1(\mathbb{R}_+)$. Then, for each
$\e>0$, the Cauchy problem \eqref{1.fad-vis} and
\eqref{1.in-fad-vis} has a unique solution $u^{\e}$ defined globally
with a uniform $BV$ bound.
Moreover, as $\e\to 0$, $u^{\e}$ converges in $L^1_{loc}$ to an
entropy solution $u\in\, BV$ to
\eqref{1.fad}--\eqref{1.in-fad}, which satisfies
\begin{align}
&\|u\|_{L^\infty(\mathbb{R}_+^2)}\leq
\|u_0\|_{L^\infty(\mathbb{R})}, \label{uniform-estimate-2}\\
&TV\{u(t)\}+ \int_0^t r(t-\tau)TV\{u(\tau)\}\,d\tau \leq C\,L, \label{bv-1} \\
&\|u(t)-u(s)\|_{L^1(\mathbb{R})}\le C |t-s|, \label{Lt-stability}
\end{align}
where $L=1+\|r\|_{L^1(\mathbb{R}_+)}$, and
$C=C(TV\{u_0\},\|u_0\|_{L^\infty})$ is a positive constant
independent of $r(t)$.
\end{thm}

Furthermore, we have
\begin{thm}\label{thm unq}
Let the resolvent kernel $r(t)$ associated with $k$ be a nonnegative
and nonincreasing function in $L^1(\mathbb{R}_+)$. Let $u, v\in
BV(\mathbb{R}_+^2)$ be entropy solutions to \eqref{1.fad} with
initial data $u_0,v_0\in BV(\mathbb{R})$, respectively. Then \be \|
u(t)-v(t)\|_{L^1(\mathbb{R})}+\int_0^t r(t-\tau) \|
u(\tau)-v(\tau)\|_{L^1(\mathbb{R})}\,d\tau\leq
L\,\|u_0-v_0\|_{L^1(\mathbb{R})}. \label{L-stability-2} \ee That is,
any entropy solution in $BV$ to \eqref{1.fad}--\eqref{1.in-fad} is
unique and stable in $L^1$. As a consequence, if $u_0$ is only in
$L^\infty$, not necessarily in $BV(\mathbb{R})$, there exists a
global entropy solution $u \in L^\infty$ to
\eqref{1.fad}--\eqref{1.in-fad}.
\end{thm}

\medskip Having established the above results, we then analyze the
case when the kernel in the scalar equation \eqref{1.fad} is a
relaxation kernel $k_\de$ that depends on a small parameter $\de>0$
so that $k_{\de}(t) \rightharpoonup (\alpha-1)\,\delta(t)$ weakly as
measures when $\de\to 0+$, where $\delta(t)$ denotes the Dirac mass
centered at the origin.
That is,
\begin{eqnarray}
&&\label{fad}
u^{\de}_t+f(u^{\de})_x+\int_0^tk_{\de}(t-\tau)f(u^{\de}(\tau))_x\,d\tau=0,\\
&&\label{in-fad} u^{\de}(0, x)=u_0(x)\in BV(\mathbb{R})
\end{eqnarray}
with $\sup_{\de>0}\|k_\de\|_{L^1(\mathbb{R}_+)}<\infty$.

We denote the entropy solution to the above problem by $u^\de(t,x)$.
Let $r_\de$ be the resolvent kernel associated with $k_\de$ via
\eqref{1.resolvent relation}. Hence, \eqref{fad} reduces to
\be\label{res-fad} u^\de_t+f(u^\de)_x+r_\de(0)u^\de=r_\de(t)
u_0-\int_0^t r_\de'(t-\tau) u(\tau)\,d\tau\ee with $
\sup_{\de>0}\|r_\de\|_{L^1(\mathbb{R}_+)}<\infty. $

By Theorems \ref{thm1}--\ref{thm unq}, we conclude that the unique
entropy solution sequence $u^\de\in BV$ to
\eqref{fad}--\eqref{in-fad}
is uniformly bounded in $L^\infty$ and uniformly $L^1$--stable,
independent of $\de$.
Then the solution sequence $\{u^\de\}$ is a compact set in
$L^1_{loc}$ so that we can extract a subsequence $\{u^{\de_k}\}$
that converges in $L^1_{loc}$ to an admissible weak solution of the
local conservation law: \be\label{cons}u_t+\alpha f(u)_x=0,\ee with
Cauchy data $u_0\in BV$.

\begin{thm}\label{thm2}
Consider the Cauchy problem \eqref{fad}--\eqref{in-fad} with
$u_0\in BV(\mathbb{R})$. Let the resolvent kernel $r_\de$ associated
with $k_\de$ as defined in \eqref{1.resolvent relation} be a
nonnegative, nonincreasing function with uniform $L^1$-norm
independent of $\de$.
Then the entropy solutions $u^\de$ to \eqref{fad}--\eqref{in-fad}
are uniformly bounded in $L^\infty$ and stable in $L^1$:
\begin{eqnarray*}
&& \|u^\de(t)\|_{L^\infty(\mathbb{R}_+^2)}\leq \|u_0\|_{L^\infty(\mathbb{R})}, \\
&&TV\{u^\de(t)\}\le C\,L,\\
&&\|u^\de(t)-u^\de(s)\|_{L^1(\mathbb{R})}\le C\, |t-s|,\\
&& \|u^\de(t)-v^\de(t)\|_{L^1(\mathbb{R})}\leq
L\|u_0-v_0\|_{L^1(\mathbb{R})},
\end{eqnarray*}
where $L:=1+\sup_{\de>0}\|r_\de\|_{L^1}<\infty$,
$C=C(TV\{u_0\},\|u_0\|_{L^\infty})>0$ is independent of $\de$, and
$v^\de(t,x)$ is the entropy solution to \eqref{fad}--\eqref{in-fad}
with initial data $v_0\in BV$.
Furthermore, if $k_{\de}(t)\rightharpoonup (\alpha-1)\, \delta(t)$
weakly as measures when $\de\to 0+$, then $u^{\de}$ converges in
$L^1_{loc}$ to an admissible weak solution $u$ of the Cauchy problem
\eqref{cons} and \eqref{1.in-fad} with initial data $u_0\in BV$.
\end{thm}

In Theorems \ref{thm1}--\ref{thm2}, the assumptions on $r_\de(t)$,
or $r(t)$, can easily be converted to the assumptions on $k_\de(t)$,
or $k(t)$, because of their symmetry between the kernel and the
resolvent through \eqref{1.resolvent relation}. For example, such
kernels $k_\de(t)$ especially include the following Set of Kernels
(i)--(vi):

(i) $k'_\de(t)\ge 0$ and $\|k_\de\|_{L^1(\mathbb{R}_+)}\le K$ for
some constant $K$ independent of $\de>0$;

(ii) $|1+\hat{k}_\de(z)|\ne 0$ for any $z$ with $Re(z)\ge 0$, and
$\hat{k}_\de(z)(1+\hat{k}_\de(z))\le 0$ for  the Laplace transform
$\hat{k}_\de$ of $k_\de$;

(iii) $sup_{\omega \in \mathbb{R}}|(1+\hat{k}_\de (i\omega))^{-1}|
\le q$ for some constant $q$ independent of $\de$;

(vi)  There exist positive numbers $T\sim \de$ and $\tau\sim \de$
such that
\[\int_{|s|\geq T} |k_\de(t)|\leq\frac{1}{12 q},
\qquad
\sup\limits_{0<s<\eta}\int_{\mathbb{R}}|k_\de(t)-k_\de(t-s)|\,dt\leq
\frac{1}{4}.\]

\smallskip The prototype is
\begin{equation}\label{exam1}
 k_\de(t)=-\frac{1-\alpha}{\de} exp(-\frac{t}{\de}), \qquad
0<\alpha<1,
\end{equation}
for which the corresponding family of resolvent kernels is
\begin{equation}\label{exam2}
r_\de(t)=\frac{1-\alpha}{\de} exp(-\frac{\alpha t}{\de}).
\end{equation}

\smallskip
In Section 2, we develop techniques for the nonlocal case, motivated
by Vol'pert-Kruzkov's techniques \cite{Volpert,K} and the
$L^\infty$-estimate techniques for the local case, to establish uniform
$L^\infty$ and $BV$ estimates of the vanishing
viscosity approximate solutions, independent of $\e$, by using the
damping nature of the memory term.

As a corollary of these estimates, we establish the convergence of
the vanishing viscosity approximate solutions to obtain the
existence of entropy solutions in $BV$. In Section 3, we show that
the entropy solution in $BV$ is unique and stable in $L^1$ with
respect to the initial perturbation. In Section 4, we prove Theorem
\ref{thm2} and discuss the hypotheses of the theorems.
Finally we give an example and show the relation of the fading
memory limit with the zero relaxation limit as first considered
systematically in Chen-Levermore-Liu \cite{CLL}; also see
\cite{LN,STW,Yo} for the model.

%%%%%%%%%%%%%%%%%%%%%%%%%%%%%
\section{Proof of Theorem \ref{thm1}}

In this section, we establish the uniform  $L^\infty$ and $BV$
estimates, as well as the uniformly continuous dependence on time in
$L^1$, which are used not only for the global existence of the
vanishing viscosity approximate solutions, but also for their
compactness. We also establish the existence and regularity of
entropy solutions.

\smallskip
\subsection{$L^\infty$ Estimate.} We first obtain a uniform $L^\infty$
estimate. Note that, by employing the resolvent kernel $r$, equation
\eqref{1.fad-vis} can be written in the form \eqref{1.res-fad-vis},
i.e. we study the Cauchy problem
\begin{eqnarray}
&& u_t^{\e}+f(u^{\e})_x+r(0)u^{\e} =r(t)
u_0^\e-\int_0^t r'(t-\tau) u^{\e}(\tau)\,d\tau+\e\, u^{\e}_{xx},\label{2.res-fad-2}\\
&& u(0,x)=u^\nu_0(x),\hspace{1cm} x\in\mathbb{R}, \label{2.res-in-2}
\end{eqnarray}
where $u^\nu_0(x)$ is defined in \eqref{1.in-z}.
By rescaling the coordinates, $(t,x)\to (s,y)=(r_0 t, r_0 x)$, we
rewrite \eqref{1.res-fad-vis} as \be\label{2.4.res-fad-vis}
\bar{u}_s^{\e}+f(\bar{u}^{\e})_y+\bar{u}^{\e}
=\frac{1}{r_0}r(\frac{s}{r_0}) u_0^\nu-\frac{1}{r_0^2}\int_0^s
r'(\frac{s-\tau}{r_0}) \bar{u}^{\e}(\tau)\,d\tau+\e\,
r_0\,\bar{u}^{\e}_{yy},\ee where $r_0>r(0)$ is any positive constant
and $\bar{u}^\e(s,y):=u^\e(\frac{s}{r_0},\frac{y}{r_0})$. For any
even integer $p$, multiplying \eqref{2.4.res-fad-vis} by $p\,
|\bar{u}^\e|^{p-1}$ and integrating over $[0,S]\times\mathbb{R}$, we
obtain
\begin{align}
\int_{-\infty}^\infty |\bar{u}^\e&(S)|^p\,dy
+ p\int_0^S\int_{-\infty}^\infty |\bar{u}^\e(s,y)|^p\,dyds\notag\\
\leq & \int_{-\infty}^\infty |\bar{u}^\e_0(y)|^p\,dy
+p\int_0^S\int_{-\infty}^\infty
\frac{1}{r_0}r(\frac{s}{r_0})|\bar{u}^\e_0(y)|
|\bar{u}^\e(s,y)|^{p-1}\,dyds
\notag\\
&-p\int_0^S\int_{-\infty}^\infty |\bar{u}^\e(s,y)|^{p-1}\int_0^s
\frac{1}{r_0^2} r'(\frac{s-\tau}{r_0})|\bar{u}^\e(\tau,y)|\,d\tau
dyds. \label{p-inq}
\end{align}
By employing the standard inequality $ab\leq \de_0
\,\frac{a^p}{p}+\de_1 \,\frac{b^q}{q}$ with
$\de_0=(\frac{4(p-1)}{p})^{p-1}$, $\de_1=\frac{p}{4(p-1)}$ and
$\frac{1}{p}+\frac{1}{q}=1$, the second term on the right-hand side
of \eqref{p-inq} is estimated as
\begin{align}
&p\int_0^S\int_{-\infty}^\infty
\frac{1}{r_0}r(\frac{s}{r_0})\,|\bar{u}^\e_0(y)|\,
|\bar{u}^\e(s,y)|^{p-1}\,dyds
\notag\\
&\quad \leq\int_0^S\int_{-\infty}^\infty \de_0
\left(\frac{1}{r_0}r(\frac{s}{r_0})\right)^p
|\bar{u}^\e_0(y)|^p\,dyds
 +\frac{1}{4} p\int_0^S\int_{-\infty}^\infty |\bar{u}^\e(s,y)|^p\,dy\,ds,
 \label{2.21}
 \end{align}
so that the last term in \eqref{2.21} is dominated by the damping in
\eqref{p-inq}. Similarly, we treat the last term in \eqref{p-inq} to
obtain
\begin{align}
\int_{-\infty}^\infty |\bar{u}^\e(S,y)|^p \,dy \leq
&\int_{-\infty}^\infty |\bar{u}^\e_0(y)|^p\,dy + \de_0\int_0^S
\left(\frac{1}{r_0}r(\frac{s}{r_0})\right)^p\,ds\,
\int_{-\infty}^\infty |\bar{u}_0^\e(y)|^p\,dy \notag\\
&+ \beta(S)\,\int_0^S\int_{-\infty}^\infty
|\bar{u}^\e(\tau,y)|^p\,dy\,d\tau,\notag
\end{align}
where
$\beta(S):=\de_0\int_0^S\left(\int_0^s\left|\frac{1}{r_0^2}r'(\frac{s-\tau}{r_0})\right|^q
\,d\tau\right)^{p/q}ds$.
By Gronwall's inequality, we have \be\label{2.4 Y}
\int_0^S\int_{-\infty}^\infty |\bar{u}^\e(\tau)|^p\,dy\,d\tau\leq
e^{ \int_0^S\beta(w)\,dw}\int_0^S\,W(s)\,ds\ee for
$W(s):=\left(1+\de_0\int_0^s
\left(\frac{1}{r_0}r(\frac{\tilde{s}}{r_0})\right)^p\,d\tilde{s}\right)\int
|\bar{u}^\e_0(y)|^p\,dy.$
 Let \be K_S=\sup\left\{
\frac{1}{r_0^2}|r'(\frac{\tilde{s}}{r_0})|:
\,\tilde{s}\in[0,S]\right\}. \ee Then, for all $s\in[0,S]$, we have
\be \beta(S)\leq \de_0K_S\int_0^S\left(-\frac{1}{r_0^2}\int_0^s
r'(\frac{s-\tau}{r_0})\,d\tau\right)^{p/q}ds \leq \de_0 K_S S,
 \ee
since $r$ is nonincreasing, and hence
$\displaystyle\lim_{p\to\infty} \frac{1}{p}\int_0^S \beta(w)\,dw=0$.
Also,
\[ \left(\int_0^S\,W(s)\,ds\right)^{\frac{1}{p}}
\leq S^{\frac{1}{p}} \left(1+\de_0^{\frac{1}{p}}
 \left(\int_0^S \left(\frac{1}{r_0}r(\frac{\tilde{s}}{r_0})\right)^p d\tilde{s}
 \right)^{\frac{1}{p}}\right)\left(\int |\bar{u}^\e_0|^p\,dy\right)^{\frac{1}{p}}.\]
Thus, if we raise \eqref{2.4 Y} to $1/p$, take the limit as
$p\to\infty$, and note that $\de_0(p)^{\frac{1}{p}}\to 4$ as
$p\to\infty$, we conclude that $
\|\bar{u}^\e\|_{L^\infty([0,S]\times\mathbb{R})}
\leq\|(1+\frac{4}{r_0}r(\frac{s}{r_0}))
|\bar{u}_0(y)|\|_{L^\infty([0,S]\times\mathbb{R})}. $ That is, \be
\|\bar{u}^\e\|_{L^\infty(\R_+^2)} \leq
\big(1+\frac{4}{r_0}\|r\|_{L^\infty(\R_+^2)}\big)
\|\bar{u}_0\|_{L^\infty(\mathbb{R})}. \ee Since $r(t)\ge 0$ is
nonincreasing so that $0\le r(t)\le r(0)<\infty$ and $r_0>r(0)$ is
an arbitrary constant, then, as $r_0\to\infty$, we conclude
$\|\bar{u}^\e\|_{L^\infty(\mathbb{R}_+^2)}\leq
\|u_0\|_{L^\infty(\mathbb{R})}$. By rescaling the coordinates
backwards, we obtain the uniform $L^\infty$--bound on $u^\e$ to
\eqref{1.res-fad-vis} independent of $\e$;
\be\|u^\e\|_{L^\infty(\mathbb{R}_+^2)}\leq
\|u_0\|_{L^\infty(\mathbb{R})}. \label{uniform-estimate-1}\ee

\smallskip
\subsection{$BV$--Regularity and Estimates}\label{S: BV}
First note that $u^\nu_0\in C^{\infty}$ has
the following bounds:
\begin{equation}
\label{2.in-bounds}\|(u^\nu_0)_x\|_{L^1}\leq
M(u_0),\qquad \|(u^\nu_0)_{xx}\|_{L^1}\le \frac{C_0}{\e} M(u_0),
\end{equation}
for some $C_0>0$ independent of $\nu$, where
$M(u_0):=TV\{u_0\}+2\|u_0\|_{L^\infty}$.

Set $v=u^\e_x$ and $w=u^\e_t$. Then the evolution equations of $v$
and $w$ are
\begin{eqnarray}&&v_t+(f'(u^\e)v)_x+\, r(0) v
= r(t) (u^\nu_0)_x-\int_0^t r'(t-\tau) v(\tau)\,d\tau+\e
v_{xx},\qquad\label{v}\\
&&v(0,x)=(u^\nu_0)_x, \label{v-in}
\end{eqnarray}
and
\begin{eqnarray}
&&w_t+(f'(u^\e)w)_x+ r(0) w =-\int_0^t r'(t-\tau)
w(\tau)\,d\tau+ \e w_{xx}, \quad\label{w}\\
&&w(0,x)=u^\e_t(0,x)
  =(u^\nu_0)_{xx}-f'(u^\nu_0)\,(u^\nu_0)_x.
 \label{w-in}
\end{eqnarray}

Multiplying \eqref{v} by $sgn (v(t,x))$ and integrating with respect
to $x$, we get
\[\dfrac{d}{dt}\|v(t)\|_{L^1}+ r(0) \|v(t)\|_{L^1}
\leq r( t) M(u_0)-\int_0^t r'(t-\tau) \|v(\tau)\|_{L^1}\,d\tau\]
since $r(0)>0$ and $r$ is nonincreasing. Integrating over
$t\in[0,T]$ yields
\begin{align}
&\|v(T)\|_{L^1}+r(0)\int_0^T \|v(t)\|_{L^1}\,dt\notag\\
&\leq \big(1+\int_0^T r( t)\,dt\big) M(u_0)-\int_0^T \int_0^t r'(t-\tau)
\|v(\tau)\|_{L^1}\,d\tau\,dt.
\end{align}
Changing the order of
integration in the last term, we arrive at
\begin{align*}
&\|v(T)\|_{L^1}+ r(0)\int_0^T \|v(t)\|_{L^1}\,dt\\
&\leq \big(1+\int_0^T r( t)\,dt\big) M(u_0) -\int_0^T
(r(T-\tau)-r(0)) \|v(\tau)\|_{L^1}\,d\tau.\end{align*} Thus, we have
\[\|v(T)\|_{L^1}+\int_0^T r (T-\tau)\|v(\tau)\|_{L^1}\, d\tau
\leq (1+\int_0^T r( \tau)\,d\tau) M(u_0).\] Because $r(\cdot)$ is
bounded in $L^1(\mathbb{R}_+)$, we obtain the following uniform
bound on the gradient $v=u_x$, \be \|u^\e_x(t)\|_{L^1}+\int_0^t r
(t-\tau)\|u^\e_x(\tau)\|_{L^1}\, d\tau\leq L\,M(u_0),
\label{bv-estimate}\ee where $L:=1+\| r\|_{L^1(\mathbb{R}_+)}$.
Similarly, we have
\be\label{2.w est} \|w(t)\|_{L^1}+\int_0^t r
(t-\tau)\|w(\tau)\|_{L^1}\, d\tau \leq \|w(0)\|_{L^1}.
\ee
Using
\eqref{2.res-fad-2} and  the bounds in \eqref{2.in-bounds} for the
initial data, we find from \eqref{w-in} that
\begin{eqnarray*}
\|w(0)\|_{L^1}\leq\nu\|(u^\nu_0)_{xx}\|_{L^1}+
\|f'(u^\nu_0)\|_{L^\infty} M(u_0)\le C_1\, M(u_0)
\end{eqnarray*}
for some $C_1>0$ independent of $\nu$. Hence, by \eqref{2.w est}, $
u^\e_t$ is uniformly bounded in $L^1$. Thus, for $0<s<t$, we get the
uniformly continuous dependence on time for the solutions to
\eqref{2.res-fad-2}--\eqref{2.res-in-2}: \be\label{2.cont dep}
\|u^\e(t)-u^\e(s)\|_{L^1}\leq \int_s^t\|w(\tau)\|\,d\tau\leq
C\,|t-s|
\ee
with $C=\max(C_1,1)\, M(u_0)$ independent of $\nu$.

\smallskip
\subsection{Existence of Entropy Solutions in $BV$
to \eqref{1.fad}--\eqref{1.in-fad}} Using
\eqref{uniform-estimate-1}, \eqref{bv-estimate}, and \eqref{2.cont
dep}, Helly's Compactness Theorem yields that a convergent
subsequence $\left\{u^{\e_m}\right\}$ may be extracted with
$\e_m\downarrow 0$ as $m\rightarrow\infty$, whose limit is denoted
by $u$, i.e., \be u^{\e_m}(t)\longrightarrow
u(t)\hspace{1cm}\text{in}\,\,L^1_{loc} \quad \mbox{for all $t>0$}.
\ee The limit $u(t,\cdot)$ is a BV function satisfying
\eqref{uniform-estimate-2}--\eqref{Lt-stability} for all $t,\,s>0$.
By construction, it is easy to check that the limit function
$u(t,x)$ is an entropy solution to \eqref{1.fad}--\eqref{1.in-fad}.

%%%%%%%%%%%%%%%%%%%%%%%%%%%%%%%%%%%%%%%%%%%%%%%%%%
\section{Proof of Theorem 1.2}
In this section, we prove the uniqueness of entropy solutions  in
$BV$ as stated in Theorem \ref{thm unq}. For any $u\in BV$, the
whole space $\mathbb{R}_+^2$ can be decomposed into three parts (see
\cite{Da,Fed,Volpert}):
$$
\mathbb{R}^2_+=J(u)\cup\mathcal{C}(u)\cup \mathcal{I}(u),
$$
where $J(u)$ is the set of points of approximate jump discontinuity,
$C(u)$ the set of points of approximate continuity of $u$, and
$I(u)$ is the set of irregular points of $u$ whose one-dimensional
Hausdorff measure is zero.

First, the entropy inequality \eqref{entropy-ineq} implies that, on
a shock $x=x(t)$ in $J(u)$, \be \sigma[\eta(u)]-[q(u)]\ge 0,\ee
where $[\eta(u)]=\eta(u(t, x(t)+0))-\eta(u(t, x(t)-0))$
and $\sigma=x'(t)$ is the shock speed.

Now assume that $u,v\in BV(\mathbb{R}_+^2)$ are the entropy
solutions with initial data $u_0, v_0\in BV(\mathbb{R})$,
respectively. Then it can be easily checked that, on $J(u)\cup
J(v)$,
$$
\sigma [ |u-v|]-[sign(u-v) (f(u)-f(v))]\le 0.
$$
In the continuous region $C(u)\cap C(v)$, since $r'(t)\leq 0$,
\begin{align*}
&\mu(t,x):=|u(t)-v(t)|_t+ q(u,v)_x+ r(0)\,|u(t)-v(t)|\\
&\qquad\qquad-r(t)|u_0-v_0|+\int_0^tr'(t-\tau)|u(\tau)-v(\tau)|\,d\tau\\
&=-\int_0^t|r'(t-\tau)|\left( |u(\tau)-v(\tau)| -sgn(u(t)-v(t))
(u(\tau)-v(\tau))\right)\,d\tau \leq 0.
\end{align*}
Therefore, $\mu$ as a measure on $\mathbb{R}_+^2$ satisfies
\begin{align*}
\mu(\mathbb{R}_+^2) =-\sum_{J(u)\cup J(v)} (\sigma[\eta]-[q])
+\mu(C(u)\cap C(v))\le 0.
\end{align*}
Then we follow the same steps as for the $BV$ estimates in Section
2.2 to conclude \eqref{L-stability-2}. When $u_0\in L^\infty$, let
$u_0^k$ be a sequence of initial data in $BV$ for which $u_0^k\to
u_0$ as $k\to\infty$. Then the $L^1$-stability result
\eqref{L-stability-2} implies that the corresponding entropy
solution sequence $u^k\in BV$ to \eqref{1.fad} with data $u_0^k$ is
a Cauchy sequence in $L^1$ which yields a subsequence converging to
$u(t,x)\in L^\infty$. It is easy to check that the limit $u(t,x)$ is
an entropy solution.

\section{Proof of Theorem \ref{thm2}}\label{S: thm2}
Let $u^\de\in BV$ denote the unique entropy solution to
\eqref{res-fad} with initial data $u_0\in BV$. Then the solution
sequence $\{u^\de\}$ is uniformly bounded and is uniformly stable in
$L^1$ with respect to the initial data since \begin{align*}
\|u^\de\|_{L^\infty(\mathbb{R}_+^2)}\le
\|u_0\|_{L^\infty(\mathbb{R})}, \qquad
\|u^\de(t)-v^\de(t)\|_{L^1}\leq L\,\|u_0-v_0\|_{L^1}
\end{align*}
and satisfies the following apriori uniform bounds: \be
TV\{u^\de(t)\}\leq C\,L, \qquad
\|u^{\de}(t)-u^{\de}(s)\|_{L^1(\mathbb{R})}\le C\, |t-s|. \ee This
implies that there exists a convergent subsequence
$\left\{u^{\de_m}\right\}$  with $\de_m\to 0$ as
$m\rightarrow\infty$, whose limit is denoted by $u$, i.e., $
u^{\de_m}(t,x)\to u(t,x)$ in $L^1_{loc}$. Then, since
$k_\de(t)\rightharpoonup (\alpha-1)\,\delta(t)$ weakly as measures
when $\de\to 0$, we conclude that $u$ is an admissible weak solution
of the Cauchy problem \eqref{cons} and \eqref{1.in-fad}.
%\be\label{3-cons}u_t+ \alpha\,f(u)_x=0\ee \be
%\label{3-in-cons} u(0,x)=u_0(x).\ee
The proof of Theorem \ref{thm2}
is complete.

\smallskip Finally we discuss some families of kernels
$\{k_\de\}$ that satisfy the assumptions stated in Theorem
\ref{thm2}.

Suppose that $k_\de\in L^1(\mathbb{R}_+)$ for all $\de>0$. Then, by
the Paley-Wiener Theorem \cite{PW}, the resolvent $r_\de$ of $k_\de$
is in $L^1(\mathbb{R}_+)$ if and only if
$$
|1+\hat{k}_\de(z)|\neq 0 \qquad\mbox{ for all } Re(z)\geq 0
$$
for the Laplace transform $\hat{k}_\de$ of $k_\de$. By extending
$k_\de$ as zero to the negative real axis, we choose the number $q$
such that
\[q\ge sup_{\omega \in \mathbb{R} }|(1+\hat{k}_\de
(i\omega))^{-1}|,\] and choose positive numbers $T$ and $\eta$
satisfying
\[\int_{|s|\geq T} |k_\de(t)| dt \leq\frac{1}{12 q},
\qquad
\sup\limits_{0<s<\eta}\int_{-\infty}^{\infty}|k_\de(t)-k_\de(t-s)|\,dt
\leq \frac{1}{4}.\] Then \be \label{res L1 bound}
\|r_\de\|_{L^1(\mathbb{R}_+)} \leq \left(8 \lceil 6q
T\|k_\de\|_{L^1(\mathbb{R}_+)}
\rceil\lceil8\|k_\de\|_{L^1(\mathbb{R}_+)}/\eta\rceil+6\right) q
\|k_\de\|_{L^1(\mathbb{R}_+)}, \ee where $\lceil s\rceil$ denotes
the smallest integer $\geq s$.

With this, for each $\de$, we can take $k_\de\in L^1$ such that
$|1+\hat{k}_\de(z)|\neq 0$ for all $Re(z)\geq 0$ and the numbers
defined above to be: $q$ independent of $\de$, $T\sim \de$ and
$\eta\sim \de$. Then, by \eqref{res L1 bound}, $r_\de$ is uniformly
bounded in $L^1(\mathbb{R})$. Furthermore, any kernel $k_\de(t)$ in
the Set of Kernels (i)--(vi) satisfies the assumptions in Theorem
1.2.

A prototype is  the family of kernels $k_\de(t)$ in \eqref{exam1}
that satisfies these assumptions. Then the corresponding family of
resolvent kernels is $r_\de(t)$ in
\eqref{exam2}
%=\frac{1-\alpha}{\de} exp(-\frac{\alpha t}{\de})$,
which fulfills the assumptions of Theorem \ref{thm2} when
$0<\alpha<1$. For this example, the scalar nonlocal equation
\eqref{1.fad}:
$$u_t+f(u)_x -\frac{1-\alpha}{\de}\displaystyle\int_0^t
e^{-\frac{t-\tau}{\de}}f(u(\tau))_x\, d\tau=0
$$
can also be written as a system of two equations: \be\label{3.rel
system} \left\{\begin{array}{l}
u_t+(f(u)-v)_x=0, \\
v_t=\dfrac{(1-\alpha) f(u)-v}{\de}.
\end{array}\right.\ee Then the range of
$\alpha\in (0,1)$ is the sub-characteristic condition. Thus, the
result of Theorem \ref{thm2} applying to this special case is
equivalent to establishing the convergence of the relaxation limit
\eqref{3.rel system} as considered in \cite{CLL,LN,STW,Yo}.

\begin{remark}
In order to obtain that the resolvent $r_\de$ of $k_\de$ is
integrable and $\|r_\de\|_{L^1(\mathbb{R}_+)}\leq 20$, it also
suffices by the Shea-Wainger Theorem \cite{GLS} to require  that a
family of kernels $\{k_\de\}$ satisfies that, for each $\de>0$,
$k_\de\in L^1_{loc}(\mathbb{R}_+)$ and is nonnegative,
nonincreasing, and convex on $\mathbb{R}_+$.
\end{remark}

\smallskip \noindent {\bf Acknowledgments}. Gui-Qiang Chen's
research was supported in part by the National Science Foundation
under Grants DMS-0505473, DMS-0244473, and an Alexandre von Humboldt
 Foundation
Fellowship. Cleopatra Christoforou's research was supported in part
by the National Science Foundation under Grant DMS-0244473. Also,
the authors would like to thank Professor Constantine Dafermos for
helpful discussions and comments.


\begin{thebibliography}{99}
\bibitem{CD} G.-Q. Chen and C. M. Dafermos, Global solutions in $L^\infty$ for a
system of conservation laws of viscoelastic materials with memory,
J. Partial Diff. Eqs. {\bf 10} (1997), 369--383.

\bibitem{CLL} G.-Q. Chen, D. Levermore, and T.-P. Liu,
Hyperbolic conservation laws with stiff relaxation terms and
entropy, Comm. Pure Appl. Math. {\bf 47} (1994), 787--830.

\bibitem{D2} C.~M. Dafermos, Development of singularities in the motion
of materials with fading memory, Arch. Rational Mech. Anal. {\bf 91}
(1986), 193--205.

\bibitem{D3} C.~M. Dafermos, Solutions in $L^\infty$ for a conservation law
with memory, Analyse Math\'{e}matique et Applications,
Gauthier-Villars, Paris, 1988, 117-128.

\bibitem{Da} C.~M. Dafermos, \textit{Hyperbolic Conservation Laws in Continuum Physics},
Springer-Verlag, Berlin, 1999.

\bibitem{Fed} H. Federer,
\textit{Geometric Measure Theory}, Springer-Verlag, New York,
{1969}.

\bibitem{GLS} G. Gripenberg, S.-O. Londen and O. Staffans,
\textit{Volterra Integral and Functional Equations}, Cambridge,
Cambridge University Press, 1990.


\bibitem{K} S. Kruzkov, First-order quasilinear equations with
several space variables, Mat. Sbornik \textbf{123} (1970), 228-255;
Math. USSR Sbornik \textbf{10} (1970), 217--273 (in English).

\bibitem{LN} T. Luo and R. Natalini, $BV$
solutions and relaxation limit for a model in viscoelasticity, Proc.
Roy. Soc. Edinburgh, {\bf 128A} (1998), 775--795.


\bibitem{MN} R. Malek-Madani and A.~J. Nohel, Formation of
singularities for a conservation law with memory, SIAM J. Math.
Anal. {\bf 16} (1985), 530--540.

\bibitem{MC} R.~C. MacCamy, A model for one-dimensional, nonlinear viscoelasticity,
Quart. Appl. Math. {\bf 35} (1977), 21-33.

\bibitem{NRT} J.~A. Nohel, R.~C. Rogers, and A.~E. Tzavaras, Weak
solutions for a nonlinear system in viscoelasticity, Commun. Partial
Diff. Eqs. {\bf 13} (1988), 309--322.

\bibitem{O} O. A. Oleinik, Discontinuous solutions of non-linear differential
equations. Usp. Mat. Nauk \textbf{12} (1957), 3--73; AMS
Translations, Ser. II, \textbf{26}, 95--172 (in English).

\bibitem{PW} R. E. A. C. Paley and N. Wiener,
\textit{Fourier Transforms in the complex Domain}, Amer. MAth. Soc.,
Providence, RI, 1934.

\bibitem{RHN} M. Renardy, W. Hrusa and J. A. Nohel,
\textit{Mathematical Problems in Viscoelasticity}, Longman, New
York, 1987.

\bibitem{STW} W. Shen, A. Tveito, and R. Winther,
On the zero relaxation limit for a system modeling the motions of a
viscoelastic solid, SIAM J. Math. Anal. {\bf 30} (1999), 1115--1135.

\bibitem{Volpert} A.~I. Vol'pert,
$BV$ Space and quasilinear equations (Russian), Mat. Sb. (N.S.) {\bf
73(115)}, 1967, 255--302.

\bibitem{Yo} W.-A. Yong, A difference scheme for a stiff system of
conservation laws, Proc. Roy. Soc. Edinburgh, {\bf 128A} (1998),
1403--1414.
\end{thebibliography}
\end{document}